\documentclass[12pt]{amsart}

\def\@typesizes{%
       \or{5}{6.5}\or{6}{7.5}\or{7}{8.5}\or{8}{11}\or{9}{12}%
       \or{10}{13}
       \or{\@xipt}{14}\or{\@xiipt}{15}\or{\@xivpt}{18}%
       \or{\@xviipt}{20}\or{\@xxpt}{24}}

\usepackage{amsthm}
\usepackage{amsmath}
\usepackage{amssymb}
\usepackage{graphicx}
\usepackage{enumerate}
\usepackage{color}
\allowdisplaybreaks
\numberwithin{equation}{section}
\numberwithin{figure}{section}

\theoremstyle{plain}
\newtheorem{theorem}{ Theorem}[section]
\newtheorem{proposition}[theorem]{ Proposition}
\newtheorem{lemma}[theorem]{ Lemma}
\newtheorem{corollary}[theorem]{ Corollary}
\newtheorem{example}[theorem]{ Example}
\newtheorem{remark}[theorem]{ Remark}
\newtheorem{definition}[theorem]{ Definition}
\newtheorem{conjecture}{ Conjecture}

\def\BET{\begin{theorem}}
\def\ENT{\end{theorem}}
\def\BEP{\begin{proposition}}
\def\ENP{\end{proposition}}
\def\BEL{\begin{lemma}}
\def\ENL{\end{lemma}}
\def\BEC{\begin{corollary}}
\def\ENC{\end{corollary}}
\def\BEE{\begin{example} \rm}
\def\ENE{\end{example}}
\def\BER{\begin{remark} \rm}
\def\ENR{\end{remark}}
\def\BED{\begin{definition} \rm}
\def\END{\end{definition}}
\def\BECJ{\begin{conjecture}}
\def\ENCJ{\end{conjecture}}

\def\bea{\begin{eqnarray}}
\def\eea{\end{eqnarray}}

\def\beas{\begin{eqnarray*}}
\def\eeas{\end{eqnarray*}}

\def\beq{\begin{equation}}
\def\eeq{\end{equation}}

\def\beal{\begin{align*}}

\def\eeal{ \end{align*} }

\def\row{\nonumber \\ & & }
\def\roweq{\nonumber \\ &=& }
\def\rowleq{\nonumber \\  & \leq & }

\def\bbC{{\mathbb C}}
\def\bbD{{\mathbb D}}

\def\bbN{{\mathbb N}}

\def\bbR{{\mathbb R}}

\def\bbZ{{\mathbb Z}}

\def\cA{{\mathcal A}}

\def\cD{{\mathcal D}}

\def\cF{{\mathcal F}}

\def\cH{{\mathcal H}}

\def\cP{{\mathcal P}}

\def\sF{{\sf F}}

\begin{document}

\title{On the Bergman projection and kernel in periodic planar domains}

\author{Jari Taskinen}

\address{Department of Mathematics and Statistics,
 University of Helsinki, P.O. Box 68, 00014 Helsinki, Finland}

\begin{abstract}
We study Bergman kernels $K_\Pi$  and projections $P_\Pi$  in unbounded planar domains $\Pi$,  which are periodic in one dimension. In the case $\Pi$ 
is simply connected we write the kernel $K_\Pi$ in terms of 
a Riemann mapping  $\varphi$ related to the bounded periodic cell $\varpi$ of 
the domain $\Pi$. 
We also introduce and adapt to the Bergman space setting the Floquet transform 
technique, which is a standard tool for elliptic spectral problems in periodic 
domains. We investigate the boundedness properties of the Floquet transform 
operators in Bergman spaces and derive a general formula
connecting $P_\Pi$ to a projection on a bounded domain. We show how this 
theory can be used to reproduce the above kernel formula 
for $K_\Pi$. Finally, we consider weighted $L^p$-estimates for $P_\Pi$
in periodic domains. 
\end{abstract}

\maketitle

\section{Introduction.}
\label{sec1}

\bigskip

In this paper we consider Bergman spaces $A^p(\Pi)$ and the Bergman 
projection $P_\Pi$ on planar domains $\Pi \subset \bbC$, which have the special 
geometry of being  periodic in  one direction; the domain is obtained as 
the union of infinitely  many translated copies of the bounded periodic 
cell $\varpi \subset \bbC$. 
We aim to work out the basic concepts of the Floquet transform 
theory in the context of operator theory in Bergman spaces of periodic domains.
If the periodic domain $\Pi$ is considered just as a domain
in the two dimensional real Euclidean space, the Floquet-transform,
acting in appropriate Sobolev spaces, is a standard tool for example in  
spectral elliptic boundary value problems, see e.g. \cite{Na1}. 
Here, we will consider the mapping properties of the Floquet-transform in Bergman 
spaces. The main goal is to study the connection of the Bergman 
projection and kernel on $\Pi$ with certain kernels related with the periodic 
cell especially  in the simply connected case. 

Let us fix some notation. Given a domain $\Omega$ in the complex plane $\bbC$, we
denote by $L^2(\Omega)$ the usual Lebesgue-Hilbert space with respect to 
the (real) area measure $dA$ and by  $A^2(\Omega)$ the corresponding Bergman space, which
is the subspace consisting of analytic functions. It is a consequence of the
Cauchy integral formula that the norm topology of $A^2(\Omega)$ is stronger 
than the topology of the uniform convergence on compact subsets, and this 
implies that the Bergman space is always a closed subspace of $L^2(\Omega)$,
hence it is a Hilbert space and in particular complete. We denote by $P_\Omega$
the orthogonal projection from $L^2(\Omega)$ onto $A^2(\Omega)$. It can 
always written by the help of the Bergman kernel $K_\Omega :
\Omega \times \Omega \to \bbC$,
\bea
P_\Omega f(z) = \int\limits_\Omega K_\Omega (z,w) f(w) dA(w)
\eea
and the kernel has the properties that $K_\Omega (z, \cdot) \in L^2(\Omega)$ 
for all $z$ and $K(z,w) = \overline{ K(w,z)}$ for all $z,w \in \Omega$. 
See e.g. \cite{K1} for a proof of these assertions.

We denote the periodic cell by  $\varpi$ and assume that $\varpi \subset ]0,1[ 
\times ]-M, M[ \subset \bbR^2 \cong \bbC$  for some $M>0$ and that its  
intersection   with the axis $\{ z \, : \, \mbox{Re}\, z = \frac12 \pm \frac12 
\}$ coincides with $\{ \frac12 \pm \frac12  \} \times (a,b) =: J_{\pm } $
for some real numbers $b > a $. 
We denote the translates of $\varpi$ by $\varpi_m = \varpi + m$, where $m 
\in \bbZ  \subset \bbC$, and then define the periodic domain $\Pi$ 
as the interior of the set
\bea
\bigcup_{m \in \bbZ} {\rm cl}( \varpi_m ).  \label{1.2}
\eea
where cl$(A)$ denotes the closure of the set $A$. Moreover, it is assumed
that $\varpi$ and $\Pi$ are  Lipschitz domains such that the boundaries $\partial 
\varpi$ and $\partial \Pi$ are in addition piecewise smooth: more precisely,  
$\partial \varpi$ consists of finitely many parametrized curves  $\gamma : [0,1] \to 
\bbC$ where $\gamma$ is continuously differentiable on $[0,1]$. Note that the 
Lipschitz assumption excludes cusps both in $\varpi$ and $\Pi$. 
Consequently,  $\partial \varpi$ is a Jordan curve. This implies that $\varpi$ is a 
Carath\'eodory domain so that by an old result of Carleman, \cite{Ca}, polynomials 
form a  dense subspace of the Bergman space $A^2(\varpi)$ (see also the 
introduction in the paper \cite{Br}). One more technical, geometric assumption $(A)$ 
will be posed in Section \ref{sec0}.
Note also  that $\Pi$ is simply connected, if and only $\varpi$ is; in the case $\varpi$ is multiply connected, 
the periodic domain is of course no more finitely connected.

Our results are as follows. In Section \ref{sec0} we consider simply
connected $\Pi$, and in Proposition \ref{prop0.4} 
we present a general formula for the kernel $K_\Pi$ in terms of a Riemann 
mapping $\varphi$ related to the  periodic cell $\varpi$.  This is actually a 
direct consequence of the fact that it is easy to write the Riemann mapping
from $\Pi$ onto the strip, where the kernel has an explicit formula. 
We consider cases of domains with polygonal boundaries with explicit 
Schwartz-Christoffel-type formulas for  the Riemann mapping.

In Section \ref{sec1a} we present the general definition of the Floquet 
transform $\sF$ for functions of real variables in periodic planar domains.
Briefly, $\sF$ transforms functions defined on the periodic domain $\Pi$ into 
functions on $\varpi$ depending on the so called Floquet variable or 
quasimomentum $\eta \in [-\pi, \pi]$. 
We observe that the definition works nicely also in the case of analytic
functions, and in particular 
in Theorem \ref{prop1.2} we prove the isometry of $\sF$ between the 
space $A^2(\Pi)$ and the corresponding parametrized Bergman-type 
space $L^2(-\pi, \pi; A_\eta^2(\varpi))$ on the periodic cell.

In Section \ref{sec2} we apply the Floquet transform techniques to formulate
a general formula involving the projection $P_\Pi$ and $\eta$-dependent
Bergman-type 
projections  $P_\eta$ related to the periodic cell $\varpi$. In contrast to Section
\ref{sec0}, this consideration is made for general, not necessarily
simply connected domains. 

Section \ref{sec2b} contains some elementary
preparations for Section \ref{sec3}, where we derive a presentation for the
projections $P_\eta$. This formula is applied in Section \ref{sec5}
to reproduce the kernel formula for $P_\Pi$ in the simply connected case by using 
Floquet transform techniques.  We also give an application of the kernel formula
to weighted $L^p$-estimates. 

As for the notation used in this paper, we let $[x]$ be the integer part of the number $x \in \bbR$, i.e. the largest integer $y$ with $y \leq x$. 
We write $\bbR^+ = \{ x + iy \in \bbC \, : \, x > 0 \}$. If $a \in \bbC$ and 
$r > 0$, $B(a,r) $ denotes the open Euclidean disc with center $a$ and radius $r$. 
Given a domain $\Omega \subset \bbC $ we denote by $\Vert f \Vert_\Omega$ the norm 
of $L^2(\Omega)$, which is
the $L^2$-space with respect to the (real) area measure. 
If $H$ is a Hilbert space, its inner product is denoted by $(f | g)_H$,
where $f,g \in H$. We recall that if $\varphi$ is a conformal 
mapping from the domain $\Omega \subset \bbC$ onto the domain $\Omega'$,
then the $|\varphi'|^2$ is the Jacobian of the coordinate transform
$\varphi$ in  the area integrals.

\section{Bergman kernels and conformal mappings on periodic domains.}
\label{sec0}

If $\Omega$ and $\Omega'$ are conformally equivalent domains in 
$\bbC$ and $\varphi : \Omega'  \to \Omega$ is a conformal mapping
and if $K(z,w)$ is the Bergman kernel of the domain $\Omega'$, then
the Bergman kernel of the domain $\Omega$ can be got from the 
formula 
\bea
K_\Omega(z,w) = K\big(\varphi(z),\varphi(w) \big) \varphi'(z) 
\overline{\varphi'(w)} , \label{0.1}
\eea
see for example \cite{Be}, formula (1), or \cite{K1}, Proposition 2.7. 
We refer to the papers \cite{BC},  \cite{FZ}, \cite{G} for general constructions of 
conformal mappings onto periodic domains, but we will here construct one which
leads to a useful formula of the Bergman kernel. 

From now on we assume in this section that $\Pi$ and $\varpi$ 
are simply connected, denote $E(z):= e^{i2\pi z}$ and define the domain $\cD := \{ 
E(z) \, : \, z \in \varpi \cup J_+ \cup J_-\}$, which  as  a consequence of the 
geometric assumptions in Section \ref{sec1} is  doubly connected. Consequently, 
there exists a conformal mapping $\phi$ from $\cD$ onto the annulus 
$\cA = \{ z \, : \, 1/\rho < |z| < \rho \}$ for a uniquely determined 
$\rho > 1$ (see for example \cite{G}, 
\cite{KK}, p. 362, 
the discussion in \cite{Sy}). We also denote
\bea
\cD_| = \{ E(z)  \, : \, z \in \varpi \}
= \cD \setminus \big( \cD \cap \bbR^+ \big) \ \ \mbox{and} \ \ 
\cA_| = \phi (\cD_|).  \label{0.1r}
\eea

We next observe that the function $z \mapsto (i 2 \pi)^{-1} \log (E(z)) $ can 
be extended as an analytic function to the periodic domain $\Pi$. This is done by 
defining the branch cut of log on the annulus $\cA$ on the curve $\Gamma:=
\phi ( \cD \cap  \bbR^+) = \phi(E(J_+)) = \phi(E(J_-))$. In order to describe this 
in detail and for the sake of simplicity, we pose our 
final additional geometric assumption, which obviously could be made weaker.

\smallskip

\noindent {\it (A) We assume that  $\Gamma$ is contained, for some $0< \delta < 1$, in the 
sector $\{z \in \bbC  \, : \,  \delta < \arg \,z < \pi - \delta \}$.}

\smallskip

\BER 
It is plain that $(A)$ is satisfied, if the mapping  $\phi$ distorts the 
arguments of  $z \in \cD \cap \bbR^+$ at most by $\pm(\pi/2 - \delta)$ for some fixed 
$\delta >0$. One then achieves $(A)$ by applying a rotation of the annulus $\cA$,
and redefining $\phi$ accordingly. This holds for example, if $\Pi$ is a small 
1-periodic perturbation of a horizontal strip. However, we conjecture that 
assumption $(A)$ is completely unnecessary for the results of this paper. 
\ENR

Let us denote $\varpi_+ = \varpi \cap \{z \in \bbC
\, : \, {\rm Re} z \geq 1/2\}$ and $\varpi_- = \varpi \setminus \varpi_+$.
If $\tilde z \in \Gamma$, then it has two preimages $z_\pm \in J_\pm$ with 
$\phi(E(z_\pm)) = \tilde z$. We consider a disc $B(\tilde z,r) \subset \cA$, where
$r = r(\tilde z)$ is so small that both sets $B(\tilde z,r ) \cap 
\phi(E(\varpi_\pm))$ are connected, and we define 
$\arg \,z \in [0, 3\pi]$ such that 
\beas
& & {\arg}\, z < \pi \ \mbox{for} \ z \in B(\tilde z,r ) \cap \phi(E(\varpi_-))
\ \ \mbox{and} 
\row
{\arg}\, z >  2\pi \ \mbox{for} \ z \in B(\tilde z,r ) \cap \phi(E(\varpi_+)).
\eeas
It is clear that 
\bea
\lim\limits_{\stackrel{\scriptstyle z \to \tilde z}{z \in\phi(E(\varpi_+))} } \arg\, z = 2 \pi + \!\!\!\!\!\! \lim\limits_{\stackrel{\scriptstyle z \to \tilde z}{z \in\phi(E(\varpi_-))} } \arg\, z .  \label{0.75}
\eea
The logarithm is defined accordingly in $\phi(E(\varpi))$,

\BEL
\label{lem0.2}
The mapping
\bea
\varphi(z)= \frac{1}{i 2 \pi} \log \big( \phi( e^{i 2 \pi z} ) \big)
+ [{\rm Re}\, z]  \label{0.4}
\eea
is conformal from $\Pi$ onto the strip $S= (-\infty,\infty) \times
\big( -\frac{\log \rho }{2\pi} , \frac{\log \rho }{2\pi} \big)$. 
\ENL

Proof. On $\varpi\ni z $ we have $\varphi(z) =  \frac{1}{i 2 \pi} \log \big( \phi( e^{i 2 \pi z} ) \big)$, which is a conformal mapping, and we also have for every $m \in \bbZ$,
\beas
\varphi(\varpi_m) = \varphi(\varpi) + m \ \ \ \mbox{with} \ \ 
\varphi(\varpi_m) \cap \varphi(\varpi_n) = \emptyset, \ \mbox{if} \ m \not= n.
\eeas
It remains to observe that $\varphi$ is continuous on $\Pi$, by the explained
branch cut and the defining  formula \eqref{0.4}, and that $\varphi$ is a
surjection, since $\phi \circ E$ is a surjection onto $\cA$. \ \ $\Box$

\bigskip

Combining this lemma with formula \eqref{0.1} yields a kernel formula
for the periodic domain $\Pi$. We denote ${\rm sech}\,z = 2/(e^z + e^{-z})$, which 
is the hyperbolic secant.

\BEP
\label{prop0.4}
We have  
\bea
\ \ \ \ K_\Pi (z,w)  &=& \widetilde K_\Pi(z,w) \frac{\pi^3}{ 4  (\log \rho)^2 } 
{\rm sech}^2 \Big( \frac{ \pi^2( \varphi(z) - \overline{ \varphi( w)} ) }{
2 \log \rho}  \Big)  ,  
\label{0.6}
\eea
where 
\bea
\widetilde K_\Pi(z,w) =
e^{ i 2 \pi (z-   \varphi(z) - \overline w +   \overline{ \varphi(w) } ) } 
\phi'( e^{ i 2 \pi z} )\overline{ \phi'( e^{  i 2 \pi  w} )  },
\label{0.6a}
\eea
\ENP

Proof. A conformal mapping from the strip $\Sigma = (-\infty,\infty) 
\times (-\pi,\pi)$ 
onto the upper half-plane is given by $z \mapsto i e^{z/2}$, and the kernel of the
upper half-plane is $-1/(\pi(z- \overline w)^2)$.  By \eqref{0.1} one obtains 
\bea
K_\Sigma(z,w) = \frac{1}{16 \pi} {\rm sech}^2\big( (z- \bar w)/4\big)  \label{BS}
\eea
Hence, \eqref{0.6} follows by applying \eqref{0.1} to the conformal
mapping $z \mapsto 2 \pi^2 \varphi(z)/ \log \rho $ from $\Pi$ onto
$\Sigma$ and observing that
$\phi(e^{i 2 \pi z}) = e^{ i 2 \pi \varphi(z)}$, see \eqref{0.4}.  \ \ $\Box$

\bigskip

We will show later, among other things, how this formula also follows from
the Floquet transform theory. We now mention a concrete example, where 
the boundary of $\Pi$ consists of straight line segments: in this case,
slightly more can be said about \eqref{0.6}, namely, the inverse $\psi$ of 
the mapping $\phi$ can be described up to some point, see \cite{BC}, \cite{FZ}. 
We assume that $\partial \Pi \cap R$, where $R = \{ z\, : \,0 < {\rm 
Re}\, z < 1  \}$ is a strip, consists of two polygonal lines $P_j$, 
$j=0, 1$ with $n_j \in \bbN$ edges and vertices at the points $z_k^{(j)}$, $
k=1, \ldots, n_j$, and we have ${\rm Re} \, z_1^{(0)} = {\rm Re} \, z_1^{(1)}=1$,
${\rm Re} \, z_1^{(n_0)} = {\rm Re} \, z_1^{(n_1)}= 0 $.  We denote the interior 
angle at $z_k^{(j)}$ by $\alpha_k^{(j)}$, where
\bea
\alpha_k^{(j)} = \pi \big( \beta_k^{(j) } + 1 \big) 
\ \ \ \ \mbox{and} \ \ 
\sum_{k=1}^{n_0} \beta_k^{(0)} = \sum_{k=1}^{n_1} \beta_k^{(0)} = 0 
\ .
\eea
According to \cite{BC}, (8.14), the inverse $\psi = \phi^{-1}: \cA \to \cD$ of the 
conformal mapping  $\phi$ 
can be defined by the formula
\bea
\psi(z) = A \int\limits_1^z
\prod_{j=0}^1 \prod_{k=1}^{n_j} \bigg( P  \Big( \frac{\zeta }{a_k^{(j)}} ,
\rho \Big) \bigg)^{\beta_k^{(j)} }   \frac{d\zeta }{\zeta}
\eea
where $A \in \bbC$ is a constant,  the points $a_k^{(j)} \in \partial \cA$ are 
the images  of the vertices $z_k^{(j)}$ and $P : \cA \to \bbC$ is 
by (8.5) of \cite{BC} the mapping
\bea
P(\zeta, \rho) = (1- \zeta) \prod_{k=1}^\infty 
\big( 1 - \rho^{2k} \zeta \big)\big( 1 - \rho^{2k} \zeta^{-1} \big) .
\eea
However, as well-known, calculating the points $a_k^{(j)}$ or the
function $\phi$ explicitly is in general impossible.

\section{Floquet transform in Bergman spaces.}
\label{sec1a}

We next introduce the operator theoretic tools. 
Since the origin of the Floquet transform, denoted here by $\sF$, is 
in the field of partial differential equations for functions of real 
variables, we put aside for a moment the complex variable, identify 
$\bbC$ with the plane $\bbR^2$, and denote the variable by $(x_1,x_2) $ 
instead of $z$; the functions are still complex valued in the following.  We 
refer to  \cite{Kbook}, \cite{Ku} for general expositions, \cite{Na1} for
applications to elliptic PDE's and e.g. \cite{CNT} for the contributions of the
author.  The Floquet transform $\sf$ with respect to the variable 
$x_1$  for  functions $f \in L^2(\Pi)$ is  given  by 
\begin{equation}
\label{1.4x}
\sF f (x_1,x_2,\eta) = \frac{1}{\sqrt{2\pi}}\sum_{m\in \mathbb{Z}}
e^{-i\eta  m }f(x_1+m , x_2),
\end{equation}
where  $(x_1,x_2)\in \varpi$ and $\eta \in[-\pi,\pi)$  
is the so called Floquet  parameter, or quasimomentum. 
As is well known (see Theorem 2.2.5 in \cite{Kbook} or Theorem 4.2 in \cite{Ku}), 
the Floquet transform  establishes an isometric isomorphism (isometric linear 
homeomorphism) between  $L^2$-spaces, 
\bea 
{\sf F} : L^2(\Pi) \to L^2(-\pi,\pi; L^2(\varpi)),  \label{1.5}
\eea 
where $L^2(-\pi,\pi;X)$ is the space of vector valued, Bochner-$L^2$ integrable 
functions on $[-\pi, \pi]$ with values  in the Banach space $X$, and it is endowed 
with the norm
\bea
\Vert g \Vert_{L^2(0,2\pi; B)}=\left(\int_{-\pi}^{\pi}
\|f(\eta)\|_X d \eta\right)^{1/2}\,.
\label{1.7}
\eea  
In particular, the series \eqref{1.4x} converges in the space $ L^2(-\pi,\pi; 
L^2(\varpi))$ for every $f \in L^2 (\Pi)$, and thus it also converges
for example pointwise for almost every $\eta$ and  $(x_1,x_2)$.

The inverse transform is also explicitly given by the formula
\bea
\sF^{-1} g (x_1,x_2) = \frac1{\sqrt{2 \pi}} \int\limits_{-\pi}^\pi
e^{ i [x_1] \eta} g(x_1 -[x_1],x_2,\eta) d \eta .  \label{1.8x}
\eea
where $g \in L^2(-\pi,\pi; L^2(\varpi))$ and $[x] $ denotes the integer part of a 
real number $x$.

In the literature there is also another variation of the definition of the Floquet 
transform: if  $f \in L^2(\Pi)$, one can also define it as
\begin{equation}
\label{1.4}
{\cF}f (x_1,x_2, \eta) 
=\frac{1}{\sqrt{2\pi}}
\sum_{m\in \mathbb{Z}} e^{-i\eta  (x_1+ m) } f(x_1+m , x_2),
\end{equation}
where $(x_1,x_2)\in \varpi$ and $\eta \in[-\pi,\pi]$.
This differs from \eqref{1.4} only by a factor, which is a 
non-constant,  unimodular, smooth function. In this case the inverse 
transform is defined on $L^2(-\pi,\pi; L^2(\varpi))$ by 
\bea
{\cF}^{-1} g (x_1,x_2) = \frac1{\sqrt{2 \pi}} \int\limits_{-\pi}^\pi
e^{ i x_1 \eta} g (x_1 -[x_1],x_2 , \eta) d \eta ,  \label{1.8}
\eea
The isometry \eqref{1.5} holds  also for $\cF$.

\BER 
\label{rem1.0}
We list some very elementary facts, some of which  are consequences of the 
isomorphism relation \eqref{1.5} and will be needed throughout the paper. 

\noindent $(i)$ 
The space $L^2(-\pi,\pi; L^2(\varpi))$ can be identified in a canonical way
with a $L^2$-space defined on a bounded Lipschitz domain of $\bbR^3$, namely
the space $L^2\big( [-\pi, \pi] \times \varpi \big)$, see \cite{Hyt}, Prop. 1.2.24. 

\noindent $(ii)$
If a sequence $(f_n)_{n=1}^\infty$ converges in $L^2(\Pi)$, then, by
definition, the sequence $({\sf F} f_n)_{n=1}^\infty$ converges in 
$L^2(-\pi,\pi; L^2(\varpi))$, and consequently, there is a subsequence
$\big( {\sf F} f_{n_k} ( \cdot, \eta) \big)_{k=1}^\infty$, which 
converges in 
$L^2(\varpi)$ for almost all $\eta \in [ -\pi, \pi]$. See \cite{RRC},
Theorem 3.12 and its vector valued version as explained after Definition 
1.2.15. of \cite{Hyt}. Apparently, this does not
need to happen for all $\eta$; a relevant case will be considered in
Example \ref{ex1.1}.
\ENR

Although we will not need Sobolev-spaces in the sequel, we mention them
briefly in order to explain the quasiperiodic boundary condition, which is
essential later. Given a domain $\Omega \subset \bbR^2 $, let $H^1(\Omega)$ 
denote the standard Sobolev-Hilbert space of complex valued $L^2(\Omega)$-functions $u$ which have weak  partial derivatives of the first order
belonging to the space $L^2(\Omega)$. According to the Sobolev embedding theorem, 
the space $H^1(\Omega)$ embeds into the space $C( {\rm cl} (\Omega) )$
of continuous functions in the closure of $\Omega$  for
example in the case $\Omega$ is a bounded Lipschitz domain.  
In particular, the functions in $H^1(\Omega)$ 
have well-defined, continuous boundary values in $\partial \Omega$. 
The Floquet transform ${\sf F }$ is also an isomorphism from the 
Sobolev space  $H^1(\Pi)$ onto  $L^2(-\pi,\pi; H^1_\eta(\varpi))$, 
where $H^1_\eta(\varpi)$ consists of Sobolev functions $u$ with 
the quasi\-per\-iodic boundary conditions
\bea 
\label{1.9}
u (1,x_2 )= e^{i\eta}u(0,x_2) , &\quad& a < x_2 < b 
\eea 
and $L^2(-\pi,\pi; H^1_\eta(\varpi))$ is defined as the closed
subspace of $L^2(-\pi,\pi; H^1(\varpi))$ such that
$f(\cdot , \eta) \in  H^1_\eta(\varpi)$ for a.e. $\eta$. 
(If the other version $\cF$ of the Floquet transform were used,
one would replace $H^1_\eta(\varpi)$  by $H_{\rm per}^1(\varpi))$,
where instead of \eqref{1.9} the periodic boundary conditions, i.e.
\eqref{1.9} with $\eta = 0$, hold.)

The terminology for the Floquet transform varies from source to source: it may
be called the Floquet-Bloch, Bloch or Zak  transform. In the Russian school 
of function spaces it is also known as the Gelfand transform, which 
is not to be mixed with the Gelfand transform of commutative 
Banach algebra theory.   In the theory of spectral elliptic boundary problems, 
the transform $\sF$ is  used
to convert a problem in the domain $\Pi$ into another one in the cell $\varpi$,
containing quasiperiodic boundary conditions. 
Here, we instead wish to study the Bergman kernel in $\Pi$ with the help of some kernels in $\varpi$.

We proceed to work in Bergman spaces of analytic functions. The transform 
\eqref{1.4}, leading to periodic boundary condition,  is not  useful as such  
for the study of analytic function spaces, since 
$\sF f$ is not analytic even if $f$ happens to be. Hence, we will use 
the transform \eqref{1.4x}, which is written using the complex variable, 
for $f \in A^2(\Pi)$,
\bea 
\sF f(z,\eta)&=& \widehat f (z ,\eta ) =  \frac{1}{\sqrt{2\pi}}\sum_{m\in \mathbb{Z}}
e^{-i \eta m  } f (z + m ), \ \ z \in \varpi, \ \eta \in [-\pi, \pi]. \label{1.14}
\eea
We note that $\sF$ is a continuous (isometric) mapping from
$ A^2(\Pi)$ into $L^2(-\pi, \pi ; L^2(\varpi))$. The series in \eqref{1.14}
converges in $ L^2(-\pi,\pi; L^2(\varpi))$, thus also pointwise
for a.e. $\eta,z$, and also in $L^2(\varpi)$ for a.e. $\eta$.  If  $g \in L^2(-\pi, \pi ; 
L^2(\varpi)$, we also denote
\bea
\sF^{-1} g (z ) & = &  \frac1{\sqrt{2 \pi}} \int\limits_{-\pi}^\pi
e^{ i [{\rm Re} z] \eta} g(z -  [{\rm Re}\,z ] , \eta ) d \eta ,  \ \ z \in \Pi  .
\label{1.14b}  
\eea 
The fact that formula \eqref{1.14b}  indeed gives the inverse of $\sF$ will be 
established only in Theorem \ref{prop1.2}, since we first need to specify
the proper function spaces. As for the motivation of the next proposition, we 
mention that the boundary condition \eqref{1.9} has to be taken into account
in the subsequent considerations even though Bergman space functions in general do 
not have well defined boundary values.

\BEP
\label{lem1.1}
$(i)$
If $f \in A^2(\Pi) \subset L^2(\Pi)$, then for almost all $\eta \in [-\pi, \pi]$,
the function $z \mapsto  \sF  f (z , \eta)$  is analytic in $\varpi \ni z $. 

\noindent
$(ii)$ There is a dense subspace $X$ of $ A^2(\Pi)$, the elements 
$f$ of which have the property that for all $\eta \in [-\pi, \pi]$,
the function $\sF  f ( \cdot, \eta)$ is analytic in $\varpi$ and has an 
analytic extension  to a neighborhood of the closure of $\varpi$ in $\Pi$.
\ENP

For an arbitrary $f \in A^2(\Pi)$, the analyticity of $\sF  f ( \cdot, \eta)$ 
does not necessarily hold for {\it all} $\eta$, see Example \ref{ex1.1} below. 

\bigskip

Proof. $(i)$ Let $0 \not=f \in A^2(\Pi)$ and $\varepsilon > 0$, and denote 
\bea
\varphi_\varepsilon (z) = e^{- \varepsilon z^2}  \label{1.35}
\eea
Then, we have $\varphi_\varepsilon f \to f$ in $A^2(\Pi)$, as $\varepsilon \to
0$. (The proof of this can be obtained by choosing a large enough compact
subset $K$ such that the $L^2$-norm of the restriction of $f$ to
$\Pi \setminus K$ is small, and then observing that $\varphi_\varepsilon$
converges to the constant function 1 uniformly on $K \cap \Pi$.) 

We now consider the domain $\varpi_U$ which is the interior of 
cl\,$( \varpi \cup \varpi_{-1} \cup \varpi_1)  \subset \Pi$.
Let us fix $z_0 \in \varpi_U$ and a small enough open 
disc $B(z_0,\varrho) $, $\varrho > 0$, such that cl\,$ 
\big( B(z_0, 3 \varrho) \big) \subset \varpi.
$ 
The Cauchy integral formula implies the estimate
$ 
|f(z)| \leq  
C \varrho^{-1} \Vert f \Vert_\varpi
$ 
for all $z \in B(z_0, \varrho)$. 
Thus, applying a translation, or the same argument in all sets 
$\varpi_U + m$, $m \in \bbZ$, we get 
\bea
\sup\limits_{z \in \varpi_U }  |f
(z + m )| = \sup\limits_{z \in \varpi_U +m} 
|f(z)| &\leq & \frac{C}{\varrho} \Vert f \Vert_{\varpi +m}
\leq 
\frac{C}{\varrho} \Vert f \Vert_{\Pi}   \label{1.40}
\eea
for all $z \in B(z_0 , \varrho)$. Moreover, 
for a fixed $\varepsilon < 1 $ we have the estimate (see the choice of the domain
$\varpi$ for the number $M$)
\bea
|\varphi_\varepsilon(z)| = e^{- \varepsilon x^2 + \varepsilon y^2}
\leq e^{M^2} e^{- \varepsilon x^2} 
\ \Rightarrow \
\sup\limits_{z \in \varpi }|\varphi_\varepsilon(z + m )|
\leq Ce^{- \varepsilon |m|^2} 
\label{1.42}
\eea
for all $m \in \bbZ$. As a consequence of these estimates, 
for a fixed
$\varepsilon$ the series in
\bea
F  ( \varphi_\varepsilon f)  (z ,\eta ) = 
\frac{1}{\sqrt{2\pi}} \sum_{m \in \bbZ} e^{- i \eta m } \varphi_\varepsilon (z + m) f(z+m)
\label{1.45}
\eea
converges uniformly in the disc $B(z_0 , \varrho)$ and thus defines an 
analytic function there. Since $z_0 \in \varpi_U$ was arbitrary, we infer
that \eqref{1.45} thus defines  an analytic function in $\varpi_U$.

Since the Floquet transform is an isometry from 
$L^2(\Pi)$ onto $L^2(-\pi,\pi;L^2(\varpi))$ and $\varphi_\varepsilon f \to f$
in $L^2(\Pi)$, we see that 
$\sF  ( \varphi_\varepsilon f)$ converges to $\sF   f$ in $L^2(-\pi,\pi;L^2(\varpi))$ as $\varepsilon \to 0$. 
This implies, in view of the definition of the vector valued norm
(see  Remark \ref{rem1.0}.$(ii)$), that there is a decreasing sequence 
$(\varepsilon_k)_{k=1}^\infty $ with $0 < \varepsilon_k \to 0$ and
\bea
\lim_{k \to \infty}  \sF  ( \varphi_{\varepsilon_k} f) (\cdot, \eta)  
= 
F   f (\cdot,\eta)  \ \ \mbox{in} \  L^2(\varpi)  \label{1.48}
\eea
for almost all $\eta \in [-\pi, \pi]$, and since 
$\sF  ( \varphi_\varepsilon f) (\cdot, \eta)$ is analytic, by \eqref{1.45}, 
the convergence
in \eqref{1.48} happens in $A^2(\varpi)$ for almost all $\eta \in [-\pi, \pi]$.
This completes the proof that $\sF  f$ is analytic in $\varpi$ 
for almost all $\eta \in [-\pi, \pi]$. 

$(ii)$ The dense subspace $X$ can be defined to consist of all functions
$\varphi_\varepsilon f$, where $f \in A^2(\Pi)$; see just below \eqref{1.35}. The result follows from
the observations around \eqref{1.45}. 
\ \ $\Box$

\bigskip
In $(i)$ of the previous proposition, the analyticity does not need
to hold for exactly all $\eta$, as shown in the next remark. 

\BEE
\label{ex1.1}
Consider the simple case that $\Pi$ is the strip $(- \infty, \infty ) \times 
(\frac14, \frac12)$, 
corresponding to $\varpi = (0,1) \times (\frac14, \frac12)$. Then, the function
$f(z) = 1/z$ is analytic on $\Pi$ and belongs to $A^2(\Pi)$. Moreover, we 
have $|z| \leq 1$ for all $z = x + iy \in \varpi$,  so that the triangle
inequality implies for $m \in \bbZ$, $m \not= 0, \pm 1$, and $z \in \varpi$,
\bea
\Big| f(z+m) - \frac{1}{m} \Big| = \Big| \frac{x+iy}{(x+iy+m)m} \Big|
\leq \frac{1}{(|m|-1)|m|} < \frac{1}{(|m|- 1)^2}. \label{1.50}
\eea
Thus for every $z \in \varpi$ and $\eta = 0$, the series 
\beas
\sum_{m \in \bbZ} e^{- i \eta m} f(z+m) = 
\sum_{m \in \bbZ} \Big( \frac{1}{m} + \Big( f(z+m) - \frac{1}{m} \Big) \Big)
\eeas
diverges by the triangle inequality and \eqref{1.50}, since 
$\sum_m 1/m$ diverges. Thus $\sF f$ is not a well defined analytic function for 
this particular value of $\eta$.
\ENE

In order to treat the Bergman projection we will need the unitarity of
certain operators, and this requires a selection of carefully 
defined function spaces. These are given in the next. 

\BED
\label{def2.0}
$1^\circ$. 
Given the number $\eta \in [-\pi, \pi]$, we denote by $A_{\eta, {\rm ext}}^2(\varpi)  $ 
the subspace of $ A^2(\varpi)$ consisting of functions  $f $ which can be 
extended as analytic functions to a  neighborhood in the domain $\Pi$ of the 
set ${\rm cl} \,(\varpi) \cap \Pi$    and satisfy  the  boundary condition 
\eqref{1.9}, i.e.
\bea
f (1 + iy ) = e^{i \eta } f ( 0 +iy )
\ \ \mbox{for all} \ a < y < b. \ 
\label{1.15}
\eea
We define the space  $A_\eta^2 (\varpi) $ as the 
closure of $A_{\eta, {\rm ext}}^2(\varpi)  $ in $A(\varpi)  $. We also denote by 
$P_\eta$ the orthogonal projection  from $L^2 (\varpi) $ onto the subspace 
$A_\eta ^2(\varpi)$.

This definition will be motivated by the isomorphy property of $\sF $ in
Theorem \ref{prop1.2}; see also the remark before \eqref{1.9} and $(ii)$ of 
Proposition \ref{lem1.1}.
Formula \eqref{1.45} gives plenty of examples of functions belonging to the 
subspace $A_{E, \eta}^2(\varpi)$. In general, the functions belonging to 
$A^2(\varpi)$  may not have properly defined boundary values on $\partial 
\varpi$ so that the  condition \eqref{1.15} cannot be posed directly.

\noindent $2^\circ$. 
We denote by  $ L^2(-\pi,\pi; A_\eta^2(\varpi))$ the subspace of $ L^2(-\pi,
\pi; A^2(\varpi))$ consisting of functions $f$ such that the function 
$z \mapsto f (z, \eta)$ belongs to $A_\eta^2(\varpi)$ for a.e. $\eta \in [-\pi, 
\pi]$. Since  $A_\eta^2(\varpi)$ is  by definition, for every $\eta$,  a closed 
subspace of $A^2(\varpi)$ and thus of $L^2(\varpi)$, we obtain that 
$ L^2(-\pi,\pi; A_\eta^2(\varpi))$ is a closed  subspace of 
$ L^2(-\pi,\pi; L^2(\varpi))$ (follows by using Remark \ref{rem1.0}.$(ii)$).

\noindent $3^\circ$. We define  the subspace 
$$
\cH_\eta := L^2 \big( -\pi,\pi;  A_{\eta,{\rm ext}}^2(\varpi) \big)
$$ 
of $L^2(-\pi,\pi; A_\eta^2(\varpi))$ which consists of functions $g$ 
such that the mapping $z \mapsto g(z ,\eta)$ belongs to 
$A_{\eta, {\rm ext}}^2(\varpi) $ for a.e. $\eta$. 
\END

It is good to keep in mind that $ L^2(- \pi,\pi; A_\eta^2(\varpi))$
and $\cH_\eta$ are not
vector valued $L^2$-spaces (since the space $ A_\eta^2(\varpi)$ depends on
$\eta$) but they only have the structure of a Banach vector bundle, 
see for example Section 1.3 of \cite{Kbook}. This fact complicates the proof of Theorem 
\ref{prop1.2}, since simple functions of $L^2 \big( -\pi,\pi;  A^2(\varpi) \big)$
are not contained in these subspaces. Thus, the following fact needs to be
proven. 

\BEL \label{lem3.9}
The subspace $\cH_\eta$ is dense in $L^2 \big( -\pi,\pi;  A_\eta^2(\varpi) \big)$
\ENL

Proof. Let us define the  Bochner space 
$L^2\big(-\pi,\pi; A_{\rm per}^2(\varpi) \big)$, where 
$ A_{\rm per}^2(\varpi) $ equals the space $A_\eta^2(\varpi)$ with $\eta = 0$, corresponding to the periodic boundary condition. 
We also denote by $ A_{\rm per, ext}^2(\varpi) \subset 
 A_{\rm per}^2(\varpi) $ the space $ A_{\eta ,{\rm per}}^2(\varpi) $ with
$\eta =0$. Finally, we define the subspace 
$L^2\big(-\pi,\pi; A_{\rm per, ext}^2(\varpi) \big)$ of
$L^2\big(-\pi,\pi; A_{\rm per}^2(\varpi) \big)$ consisting of those functions
$g$ which have values in $ A_{\rm per, ext}^2(\varpi) $ for almost every $\eta$.

By the theory of Bochner spaces, see Lemma 1.2.19 in \cite{Hyt}, a
dense subspace of $L^2\big(-\pi,\pi; A_{\rm per}^2(\varpi) \big)$ is formed
by simple functions 
\bea
\sum_j \chi_j(\eta) g_j(z),  \label{3.na}
\eea
where $\chi_j$ is the characteristic function of some measurable subset
of $[-\pi,\pi]$, $g_j \in A_{\rm per}^2(\varpi)$ and the sum contains only  finitely
many terms. Moreover, since $ A_{\rm per, ext}^2(\varpi) $ is dense in
$ A_{\rm per}^2(\varpi) $ by definition, a dense subspace of 
$L^2\big(-\pi,\pi; A_{\rm per}^2(\varpi) \big)$ is formed by functions
\eqref{3.na}, where every $g_j$ belongs to $ A_{\rm per, ext}^2(\varpi) $.
This shows that $L^2\big(-\pi,\pi; A_{\rm per, ext}^2(\varpi) \big)$ is
dense in $L^2\big(-\pi,\pi; A_{\rm per}^2(\varpi) \big)$.

The assertion of the lemma now follows from the facts that the mapping
$f(z,\eta)  \mapsto e^{i \eta z} f(z,\eta)$ is an isomorphism
from 
$L^2\big(-\pi,\pi; A_{\rm per}^2(\varpi) \big)$
onto  $L^2\big(-\pi,\pi; A_\eta^2(\varpi) \big)$
and from  $L^2\big(-\pi,\pi; A_{\rm per, ext}^2(\varpi) \big)$
onto $\cH_\eta$. The proofs of these isomorphisms are left to the reader. 
\ \ $\Box$

\bigskip

We are finally prepared to prove the main result of this section.

\BET
\label{prop1.2}
Floquet transform $\sF$ maps $A^2(\Pi)$ onto $ L^2(-\pi,\pi; 
A_\eta^2(\varpi))$. Its inverse 
$\sF^{-1} :  L^2(-\pi,\pi; A_\eta^2(\varpi)) \to A^2(\Pi)$ 
is given by the formula \eqref{1.14b}. Moreover, 
$\sF$ preserves the inner product and is thus a unitary operator. 
\ENT

Proof. We recall that ${\sf F}$ is an isometric isomorphism from $ L^2(\Pi)$ 
onto $L^2(-\pi,\pi; L^2(\varpi))$ and first show that it maps the space  
$A^2(\Pi)$ onto $ L^2(-\pi,\pi;  A_\eta^2(\varpi))$. Let the 
dense subspace $X \subset L^2(\Pi)$ be as in Proposition \ref{lem1.1}.$(ii)$
and assume  $f \in X$ so that  for every $\eta$, $\sF f (z, \eta)$ is an 
analytic function of $z$ on $\varpi$ and also has an analytic extension
to a neighborhood of the closure of $\varpi$, see Definition \ref{def2.0}.
Thus $\sF f(\cdot, \eta)$ also has well defined values on the lateral
boundaries of $\varpi$. Moreover, given 
$\eta \in [-\pi,\pi]$, the function $\sF  f ( \cdot,\eta)$ satisfies the 
quasiperiodic boundary condition \eqref{1.15}, since 
\beas 
& &  \sF  f(1 +iy ,\eta ) =
\frac{1}{\sqrt{2\pi}}\sum_{m\in \mathbb{Z}}
e^{-i \eta ( m -1)  } f (1 +iy + (m-1) )
= e^{i \eta}  \sF  f (iy ,\eta ) .
\eeas 
Hence, $\sF  $ maps the dense subspace $X$ into $L^2(-\pi,\pi; A_\eta^2(\varpi))$,
and since $\sF$ is an isometry and $L^2(-\pi,\pi; A_\eta^2(\varpi))$
is complete, it also maps $ A^2(\Pi)$ into $L^2(-\pi,\pi; A_\eta^2(\varpi))$.

To see that  $\sF $ is surjective we first consider a function $g \in  \cH_\eta$. 
Clearly, for all $m \in \bbZ$ and a.e. $\eta \in [-\pi,\pi]$, the function 
\bea
G_{\eta,m} : \varpi_m \to \bbC, \ \ G_\eta(z) =  e^{i m  \eta  }  
g( z -  m , \eta),   \label{1.60}
\eea
is analytic on the subdomain $\varpi_m$, $m\in \bbZ$ (see \eqref{1.2}).
Moreover, by Definition \ref{def2.0}, it has an analytic extension, 
still denoted by $G_{\eta,m}$, to a  neighborhood  of the closure of $\varpi_m$ in $\Pi$.   Since  $g(\cdot ,\eta)$
satisfies the boundary condition \eqref{1.15}, we get
for all $m \in \bbZ$, all $z = m +1 + iy$ with $y \in [a,b]$,
\beas
& &G_{\eta , m} (z) = \lim\limits_{\varepsilon \to 0^+}
G_{\eta,m} (m + 1 - \varepsilon +iy )
= 
\lim\limits_{\varepsilon \to 0^+} e^{i\eta m   } 
g(1 +iy - \varepsilon, \eta)
\roweq
e^{i\eta m  } g(1 +iy , \eta) = e^{i \eta( m + 1 )  }
g(iy , \eta )
= 
\lim\limits_{\varepsilon \to 0^+} e^{i\eta( m + 1 ) } 
g(iy + \varepsilon , \eta)
\roweq 
\lim\limits_{\varepsilon \to 0^+}
G_{\eta,m +1 } (m + 1  +  \varepsilon +iy )
= G_{\eta,m+1} (m +1  +iy ) = G_{\eta, m+1} (z).
\eeas
In other words, $G_{\eta,m}$ and $G_{\eta,m+1}$
are two analytic functions with overlapping domains of definition and they 
coincide on  the line  segments $J_m = \{ m +1  \} \times [a,b]$.
Thus, $G_{\eta,m}$ and $G_{\eta,m+1}$ coincide in their common domain of definition, and we can define the analytic function $G_\eta$ 
on $\Pi$ by setting
$ 
 G_\eta (z) = G_{\eta,m}(z) 
$ 
for $z$ belonging to the closure of $\varpi_m$ in $\Pi$. 

Comparing \eqref{1.60} with \eqref{1.14b} we see that 
\beas
\sF^{-1} g (z)  =  
\frac{1}{\sqrt{2\pi}}\int\limits_{- \pi}^\pi G_\eta (z) d\eta ,
\eeas
which is an analytic function in $\Pi$, since $G_\eta $ is. Hence,  $\sF ^{-1} 
g \in A^2(\Pi)$. Then, 
we have $\sF   \sF^{-1}  g = g$ for all $g \in \cH_\eta \subset L^2\big(-\pi,\pi;
L^2(\varpi)\big)$, by the general Floquet inversion formula \eqref{1.4}. 
We conclude that image of $\sF \big( A^2(\Pi) \big) $ contains the subspace $\cH_\eta$, which is dense in $L^2(-\pi,\pi; A_\eta^2(\varpi))$ by Lemma
\ref{lem3.9}. This implies the surjectivity, since ${\sf F}$ is an isometry 
and thus bounded from below with respect to the  $L^2$-norms.  

It is well known that the Floquet transform preserves the inner product. 
This can also seen directly, since for $f,g \in L^2(\Pi)$ 
and $H:= L^2(-\pi, \pi;L^2(\varpi) )$ there holds
\bea
& & ( \sF f | \sF g )_H= \sum_{m,n \in \bbZ} \int\limits_\varpi
\int\limits_{-\pi}^\pi
e^{- i m\eta} f(z+m)
e^{ i n \eta} \overline{g(z+n) } dz d\eta
\roweq
\sum_{m\in \bbZ}\int\limits_\varpi f(z+m)  \overline{g(z+m) } dz 
= \int\limits_\Pi f(z)  \overline{g(z) } dz = (f|g)_{L^2(\varpi)}. 
\ \  \Box
\eea

\section{General kernel formula for the periodic domain.}
\label{sec2}

Our purpose is to show how the Bergman projection in $\Pi$ can be presented
with the help of an orthogonal projection in the periodicity cell $\varpi$ 
by using the Floquet transform. 

\BED
\label{lem2.1} 
For all $\eta \in  [-\pi,\pi]$  we denote by $P_\eta$ the orthogonal projection from
$L^2 (\varpi)$ onto $A_{\eta }^2 (\varpi)$; see Definition
\ref{def2.0}.$1^\circ$ for the notation. 
\END

The well-known proof for the existence of Bergman kernels, see e.g. \cite{K1},
p.1060, applies also in this case and implies that $P_\eta$ can
be written with the help of  the integral kernel $K_\eta : 
\varpi \times\varpi \to \bbC$
\bea
P_\eta f (z) = \int\limits_\varpi K_\eta(z,w) 
f(w) dA(w) , \label{1.17}
\eea
where $K_\eta(z, \cdot) \in L^2(\varpi)$ for all $z \in \varpi$. 
In Section \ref{sec3} we will calculate the kernel $K_\eta$ in some cases.
We also define an operator $\cP$ by denoting, for all $f \in L^2(-\pi,\pi; 
L^2(\varpi))$,
\bea
\cP f(z,\eta) = \big( P_\eta f( \cdot , \eta ) \big)  (z) , \ \ \ z \in \varpi, \ \eta \in [-\pi, \pi].  
\eea
It follows directly from the definitions of the spaces and their norms
that this is a bounded operator from $ L^2(-\pi,\pi; L^2(\varpi))$ into  $L^2(-\pi,\pi; A_\eta^2(\varpi))$. More precisely, the following holds. 

\BEL
\label{lem2.2} The operator $\cP$ is the orthogonal projection from
$L^2(-\pi,\pi; L^2(\varpi))$  onto $L^2(-\pi,\pi; A_\eta^2(\varpi))$. 
\ENL
 
Proof. Since $P_\eta$ is the orthogonal projection  $L^2 (\varpi)$ onto $A_{\eta }^2 (\varpi)$, we have $\cP^2 = \cP$ and $\cP f = f$ for every 
$f \in L^2(-\pi,\pi; A_\eta^2(\varpi))$. Thus, $\cP$ is  a projection
operator between the given spaces, and it suffices to observe that
it is also self-adjoint. With $H:= L^2(-\pi, \pi;L^2(\varpi) )$ and 
$f,g \in H$  we obtain, by the self-adjointness of $P_\eta$,
\beas
& & (\cP f | g )_H 
=  \int\limits_\varpi \int\limits_{-\pi}^\pi
P_\eta f ( \cdot ,  \eta) \overline{g(\cdot, \eta) } d\eta dA
= \int\limits_\varpi \int\limits_{-\pi}^\pi f(\cdot , \eta)
\overline{ P_\eta g (\cdot , \eta)}  d\eta dA
= (f | \cP  g )_H.
 \ \  \Box
\eeas

We can now describe the general relation of the orthogonal  projections in the bounded domains and the full domain $\Pi$. 
 
\BET
\label{th2.3}
The Bergman projection $P_\Pi$ from $L^2(\Pi)$ onto $A^2(\Pi)$ can be written as  
\bea
& &  \sF^{-1} \cP  \sF  f (z) = 
\frac1{\sqrt{2 \pi}} \int\limits_{-\pi}^\pi e^{ i [{\rm Re} z ]\eta}
\big( P_\eta  \widehat f ( \cdot , \eta) \big) ( z -  [{\rm Re} z])  
d \eta 
\roweq
\frac1{2 \pi} \int\limits_{\Pi} \int\limits_{-\pi}^\pi 
e^{ i \eta ( [{\rm Re} z]-  [{\rm Re} w ])}
K_\eta( z -  [{\rm Re} z] , w - [{\rm Re}w ]) 
f(w) d\eta d A(w)  
\label{4.4zz}
\eea
\ENT

Proof. The operator $\sF^{-1} \cP \sF $ is self-adjoint, since
$\sF $ preserves the inner product (Theorem \ref{prop1.2}), and also clearly a projection, 
since $\cP$ is. 
So, this operator must be the orthogonal projection from $L^2(\Pi)$ onto $A^2(\Pi)$,
in view of Theorem \ref{prop1.2}. There remains to verify the rest of 
the formula \eqref{4.4zz}. 

The series \eqref{1.14} for $\sF  f$  converges 
in $L^2(\varpi)$ for almost all $\eta \in [-\pi, \pi]$, and since $P_\eta$ 
is a bounded operator, we can write
\bea
\cP  \sF  f (z;\eta)& = &  \frac{1}{\sqrt{2\pi}}\sum_{m\in \mathbb{Z}}
e^{-i \eta m } \int\limits_\varpi K_\eta(z ,w)  f (w + m )  
dA(w) , \ \ z \in \varpi,  \label{2.10}
\eea
for almost all $\eta$. On the other hand, the series \eqref{1.14} for $\sF  f$
converges also in $L^2(-\pi,\pi;L^2(\varpi))$, hence the series \eqref{2.10}
converges in this space, too. Since $\sF ^{-1}$ is a bounded operator
$L^2(-\pi,\pi;L^2(\varpi)) \to L^2(\Pi)$, the result
follows from the defining formula for $\sF^{-1}$ in \eqref{1.14b}, 
\beas
& &  \sF^{-1} \cP  \sF f(z)
\roweq  
\frac1{2 \pi} \sum_{m \in \bbZ} \int\limits_\varpi  \int\limits_{-\pi}^\pi 
e^{i \eta ( [{\rm Re} z] -m)}
K_\eta( z -  [{\rm Re} z] , w ) f(w + m) d\eta d A(w) 
\roweq  
\frac1{2 \pi} \sum_{m \in \bbZ} \int\limits_{\varpi_m} \int\limits_{-\pi}^\pi 
e^{ i \eta ( [{\rm Re} z]-  [{\rm Re} w])}
K_\eta( z -  [{\rm Re} z] , w -  [{\rm Re} w]) f(w)   d\eta d A(w)  
\roweq  
\frac1{2 \pi} \int\limits_{\Pi} \int\limits_{-\pi}^\pi 
e^{ i \eta ( [{\rm Re} z]-  [{\rm Re} w])}
K_\eta( z -  [{\rm Re} z] , w - [{\rm Re}w]) 
f(w)  d\eta d A(w)   . \ \ \Box
\eeas

We will show in Section \ref{sec3} that  if the periodic cell $\varpi$ is 
simply connected, the kernel $K_\eta$ can be constructed by applying a
certain canonical Riemann  conformal  mapping of doubly
connected domains and derive again the kernel formula of Section 
\ref{sec0} in this way.

\section{Weights, domains and projections.}
\label{sec2b}
In this section we collect some basic facts concerning the Bergman projection
in domains and weighted spaces. These results are known at least to experts in 
this area, but we give some proofs for the sake of the completeness of the 
presentation. 

Let $\Omega \subset \bbC$ be a domain. 
By  a weight on $\Omega$  we mean a continuous function $V : \Omega \to \bbR^+
=(0, \infty)$, and we denote by $L_V^2(\Omega)$
the weighted $L^2$-space on $\Omega$ with norm and inner product
\bea
\Vert f \Vert_{\Omega,V}^2 = \int\limits_\Omega |f|^2 V \, dA 
\ \ , \ \ (f |g)_{\Omega,V}=  \int\limits_\Omega f \overline g V \, dA ,  \label{1.a}
\eea
where $f,g, \in L_V^2(\Omega)$. 
We also denote by $A_V^2(\Omega)$ the subspace of $L_V^2(\Omega)$
consisting of analytic functions. The proof of the unweighted case applies
here, too, and shows that $A_V^2(\Omega)$ is a closed  subspace.
Our  weight functions will mostly be of the form 
\bea
V(z) = |v(z)|^2 = v(z)  \overline{v(z)},  \label{2.10y}
\eea 
where $v : \Omega \to \bbC \setminus \{ 0\}$  is analytic.

We will  use the following observations.

\BEL
\label{lem2.5}
Let $\Omega$ and $\Omega'$ be conformally equivalent domains and let 
$\varphi$  be a conformal mapping from $\Omega$ onto $\Omega'$. Then,
the composition operator $I: f \mapsto f \circ \varphi$ is a unitary
isomorphism from the space $L_V^2 (\Omega')$ onto $L^2 (\Omega)$
and also from $A_V^2 (\Omega')$ onto $A^2 (\Omega)$, where the weight is
$V(z) = |\psi'(z)|^2$ with $\psi = \varphi^{-1}: \Omega' \to \Omega$.

If $X$ is a closed subspace of $L_V^2(\Omega')$ and $P_X$ is the orthogonal 
projection from $L_V^2(\Omega')$ onto $X$, then 
the orthogonal projection $P$ from $L^2(\Omega)$ onto
$I (X) \subset L^2(\Omega)$ is given by $P=I P_X I^{-1}$. If the
function $K_X \in L^2( \Omega' \times \Omega' )$ is the integral kernel of
$P_X$, then the kernel of $P$ is given by
\bea
K(z,w) =  K_X \big( \varphi(z) , \varphi(w) \big) |\varphi'(w)|^2, \ \ \ \ z,w \in \Omega.
\eea
\ENL

Proof. The unitary isomorphism property of $I$ is a direct consequence of
the definitions and the fact that $V$ is the Jacobian of the conformal 
transform $\psi$. Note that the inverse of the operator $I$ is given
by $I^{-1}: L^2(\Omega) \to L_V^2(\Omega')$, $f \mapsto f \circ \psi$.

Next, it is straightforward to see that the mapping $P=IP_XI^{-1} : L^2(\Omega) 
\to  L^2(\Omega)$ is selfadjoint (since the adjoint of $I$ is $I^{-1}$ by the 
unitarity), a  projection ($P^2= P$), and maps $L^2(\Omega)$  onto $I(X)$. 
Hence, it is the claimed orthogonal projection. Finally, given $f \in 
L^2(\Omega)$ we write
\bea
Pf(z)&=&IP_XI^{-1}f(z) = \int\limits_{\Omega'} K_X \big(  \varphi (z)  , w \big) 
f \big(\varphi^{-1} (w) \big)dA(w)  
\roweq
\int\limits_{\Omega} K_X \big(  \varphi (z) , \varphi(w) \big) 
f (w) |\varphi'(w)|^2 dA(w). \hskip3cm \Box
\eea

\section{Simply connected periodic domain.}
\label{sec3}

In the following we construct the kernel for the operator $P_\eta$, assuming
eventually that $\Pi$ is simply connected and satisfied assumption $(A)$
of Section \ref{sec0}.
Due to the quasiperiodic boundary condition in the definition
of the target space $A_\eta^2(\varpi)$, the Riemann mapping from $\varpi$ onto 
$\bbD$ is not useful but instead we will use the doubly connected exponential image 
of  $\varpi$ and its Riemann map to an annulus. In the next section we will  rewrite 
the formula of Theorem \ref{th2.3} of the Bergman projection $P_\Pi$ with the
help of the mentioned Riemann map of the bounded domain and 
show how formula \eqref{0.6} follows from these considerations. 

We proceed to the construction of  the kernel $K_\eta$ of the orthogonal 
projection $P_\eta : L^2(\varpi) \to A_\eta^2(\varpi)$ for any $\eta \in [- 
\pi, \pi]$, see Definition \ref{def2.0}, $2^\circ$. This will be accomplished 
in several steps, where we also employ certain classical methods in complex 
analysis. We will use the notation of Section \ref{sec0} and recall
that the exponential map $E: z \mapsto e^{i 2 \pi z}$ maps the set $\varpi \cup J_+ \cup J_-$ onto the domain $\cD$, which 
is contained in an annulus so that there exist numbers $0 < \rho_0 < \rho_1$
with 
\bea
\cD \subset \{ z \in \bbC \, : \, \rho_0 < |z| < \rho_1 \} .  \label{3.4y}
\eea
The domain $\cD$ is doubly connected, if $\varpi$ is simply connected. 
Let $\cD_| = E(\varpi) \subset \cD$ be as in \eqref{0.1r}. We will soon employ the branch cut
of the logarithm as defined in Section \ref{sec0}, but we need to agree that 
when considered in $\cD$, the branch cut of the logarithm 
happens on  $\cD \cap \bbR^+$ so that Im\,$\log z \in (0, 2 \pi)$ for 
$z \in \cD_|$ and log is analytic on $\cD_|$. Let us denote $L = E^{-1} : \cD_| \to \varpi $, 
$L(z) = (i 2 \pi)^{-1} \log z$ . 

We define for all $\eta \in [-\pi,\pi]$ the weight function 
\bea
W (z) = |L'(z)|^2 = \frac{1}{4 \pi^2 |z|^2} \ , \ \ \ z  \in \cD,
\label{3.4}
\eea 
and consider the weighted Bergman space $A_W^2(\cD_|)$ 
defined below \eqref{1.a} (restricting the weight $W$ onto $D_|$, of course). 
The Bergman space  $A^2_{W} ( \cD)$ can be considered
as a subspace $A^2_{W}( \cD_|)$, which  consists of  functions that can 
be extended from $ \cD_|$ to the full domain $\cD$ as analytic
functions. Since $A^2_{W}( \cD)$ is complete, it is a closed subspace 
of  $A^2_{W}( \cD_|)$.

\BEL \label{lem3.}
$(i)$ The operator $I_1 : f \mapsto f \circ L$ is a unitary isomorphism 
$A^2(\varpi)  \to  A_{W}^2(\cD_|)$ and $ L^2(\varpi)   
\to L_{W}^2(\cD).$

\noindent $(ii)$ 
The operator $I_1$ maps the space $A_\eta^2(\varpi)$ 
onto $A_{W,\eta}^2(\cD_|) $, which is the closure of the subspace
\bea
A_{W,\eta,{\rm ext}}^2(\cD_|) = 
\Big\{ z^{\eta/(2 \pi)} g  \, : \, g \in A_W^2(\cD) 
\Big\} 
\label{3.12}
\eea
in $A_W^2(\cD_|)$.
\ENL

Proof. The claim $(i)$ follows by applying Lemma \ref{lem2.5} to the 
inverse operator $I_1^{-1} : f \mapsto f \circ E$.
As for the claim $(ii)$, given an arbitrary $g \in A_W^2(\cD)$, the function
\bea
f := \big(z^{\eta / (2\pi)} g \big) \circ E 
=: h \circ E \in L^2(\varpi) ,
\eea
belongs to the subspace $A_{\eta,{\rm ext}}^2(\varpi)$ of $A_\eta(\varpi)^2$ (see
Definition \ref{def2.0}) and  we have $I_1 f = h$, where  $h$ is an arbitrary 
element of $A_{W,\eta,{\rm ext}}^2$. Thus, the operator $I_1$ is surjection onto
$A_{W,\eta,{\rm ext}}^2$ and consequently also onto $A_{W,\eta}^2$, by $(i)$ 
and the closedness of $A_\eta^2 (\varpi)$ in $A^2(\varpi)$.

On the other hand, in view of Definition \ref{def2.0}, we consider a
function $f \in A_{\eta, {\rm ext}}^2(\varpi)$ and thus has well defined 
continuous  boundary values satisfying \eqref{1.15}. 
The function $f \circ L \in A_W^2(\cD_|)$ has a discontinuity at the slit
of $\cD_|$, more precisely
\bea
\lim\limits_{y \to 0^-} f \circ L (x+ iy) = e^{i \eta} 
\lim\limits_{y \to 0^+} f \circ L  (x+ iy) \ \ \mbox{for all}
\ x \in \cD \cap \bbR^+.
\eea
Hence, 
the function $g = 
z^{-\eta/(2 \pi)} f \circ L (z)$ has a continuous extension to $\cD$, which 
makes it into an  analytic function in $\cD$. Moreover, since the multiplier
$z^{-\eta/(2 \pi)}$  is a bounded function, we also have
$g \in L_W^2(\cD)$ and thus $g \in A_W^2(\cD)$, hence, we obtain $I_1 f = 
z^{\eta/(2\pi)} g \in A_{W,\eta, {\rm ext}}^2(\cD_|)$. Since, by definition,   
$A_\eta^2(\varpi)$ is the closure  of $ A_{\eta, {\rm ext}}^2(\varpi)$ in 
$A^2(\varpi)$  and  $I_1$ is an isometry, we obtain that $I_1$ maps 
$A_\eta^2(\varpi)$ into  $A_{W,\eta}^2(\cD_|) $.   \ \ $\Box$

\bigskip

From now on we will additionally assume that $\Pi$ and $\varpi$ are simply 
connected (thus $\cD$ is  doubly connected).
We recall from Section \ref{sec0} the conformal mapping
\bea
\phi: \cD \to \cA = \{ z \in \bbC \, : \, 1/\rho < |z| < \rho \}  \label{3.20}
\eea  
where  the number $\rho > 1$ is uniquely determined by $\cD$. Let us also agree 
that when considered on $\cA$, the logarithm is defined as in Section
\ref{sec0}. Note that our assumptions on the geometry of $\Pi$ and $\varpi$, 
\eqref{1.2}, do  not imply
smoothness of the boundary $\partial \cD$ of $\cD$. Hence, $\partial \cD$ may
include corners for example. This has the consequence that the derivative
$\phi'$ does not need to be bounded or bounded away from zero. See for example
\cite{P}, \cite{S1}, \cite{S2}, \cite{Be}, \cite{T3}.

We define on $\cA$ the weight 
\bea
V(z) = W(\psi (z) ) |\psi'(z)|^2 , \ \ z \in \cA,
\label{3.22}
\eea
where $\psi = \phi^{-1}:   \cA \to \cD$ is the inverse map. 
According to \eqref{3.4}, 
\bea
V(z)&  = & |(L\circ \psi)'(z)|^2 =
\frac{|\psi'(z)|^2 }{4 \pi^2| \psi(z)|^{2 } }=  v(z) \overline{v(z)}
\ \ \mbox{with}  \ \ 
v(z) = 
\frac{1}{2 \pi}\frac{\psi'(z)}{\psi(z)}  ,
\label{3.24}
\eea
where $z \in \cA .$ 
Moreover, we denote by $A_{V,\eta, {\rm ex} }^2 (\cA_|)$ the subspace of 
all functions of the form
\bea
z^{\eta/(2\pi)} g, \ \ \ \mbox{where $g \in A_V^2(\cA)$},
\label{3.27b}
\eea
and by $A_{V,\eta}^2(\cA_|)$ the closure of $A_{V,\eta, {\rm ex} }^2 (\cA_|)$ in 
$A_V^2 (\cA_|)$. We recall that the functions belonging to $A_{V,\eta}^2(\cA_|)$
may be discontinuous on the curve $\Gamma = \phi(\cD \cap \bbR^+)$, see Section
\ref{sec0}. Finally, we denote by $I_2$  the composition operator
$I_2 : f \mapsto f \circ \psi$. 

\BEL \label{lem3.4} 
$(i)$ The operator $I_2$ is a unitary isomorphism 
$L_{W}^2(\cD)  \to L_{V}^2(\cA)$. 

\noindent $(ii)$ The operator $I_2$ maps the space $A_{W,\eta}^2 (\cD_|)$
onto $A_{V,\eta}^2 (\cA_|)$ for every $\eta \in [- \pi, \pi]$. 
\ENL

As a consequence, the composition 
operator
\bea
I_2 I_1 : f \mapsto f \circ L \circ \psi
\label{3.24r}
\eea
is  a unitary isomorphism from $L^2(\varpi)$ onto $L^2_V(\cA)$
and from $A_\eta^2(\varpi)$ onto $A_{V,\eta}^2 (\cA_|)$.

\bigskip

Proof. 
Here, the first assertion is just a variant of Lemma 
\ref{lem2.5} with a similar proof. To prove $(ii)$ we note that 
due to choice of the branch cut on $\cA$ (see \eqref{0.75}) and \eqref{3.4y}, 
the function $F_\eta(z) = \big( \phi(z)/z\big)^{\eta/(2\pi)}$
is analytic, bounded and bounded away from zero on $\cD$ for all $\eta \in 
[-\pi, \pi]$, and  the same is thus true also for the function
\beas
\Big( \frac{ \psi(z)}{z} \Big)^{\eta/(2\pi)} = 
\frac{1}{F_\eta \circ \psi (z) } 
\ \ \
\mbox{on} \ \cA.
\eeas 
Given 
$f = z^{\eta/(2\pi)} g \in A_{W,\eta,{\rm ext} }^2 (\cD_|)$ with
$g \in A_W^2 (\cD)$ as in \eqref{3.12} we thus have
\bea
I_2 f(z) = \psi(z)^{\eta/(2\pi)} g(\psi(z)) =
z^{\eta/(2\pi)} \Big(\frac{\psi(z)}{z} \Big)^{\eta/(2\pi)} g(\psi(z))
=: z^{\eta/(2\pi)} h(z)
\eea
where 
\bea
& & \Vert h \Vert_{\cD,V}^2 = \int\limits_\cA \Big| \frac{\psi(z)}z 
\Big|^{\eta/(2\pi)} 
|g(\psi(z))|^2 \frac{|\psi'(z)|^2 }{4 \pi^2| \psi(z)|^{2 } } dA(z) 
\rowleq
C \int\limits_\cA |g \circ \psi|^2 
\frac{|\psi'|^2 }{4 \pi^2| \psi|^{2 } }  dA
= C \int\limits_\cD |g|^2 W  dA .
\eea
Hence, $h \in A_V^2(\cA)$, and $I_2$ maps $A_{W,\eta,{\rm ext} }^2 (\cD_|)$ into 
$A_{V,\eta,{\rm ext} }^2 (\cA_|)$. The converse relation
\beas
I_2^{-1} \big( A_{V,\eta,{\rm ext} }^2 (\cA_|) \big) \subset 
A_{W,\eta,{\rm ext} }^2 (\cD_|)
\eeas
can be proved in the same way, starting by $f = z^{\eta/(2\pi)} g\in 
A_{V,\eta,{\rm ext} }^2 (\cD_|)$ with $g \in A_V^2 (\cD)$ and using $\phi$
instead of $\psi$. 

The proof is completed by observing that the operator  $I_2 : A_{W,\eta}^2 (\cD_|) 
\to A_{V,\eta}^2 (\cA_|)$ as well as its inverse are bounded and by taking into 
account the densities of $A_{W,\eta,{\rm ext} }^2 (\cD_|)$ and
$A_{V,\eta,{\rm ext} }^2 (\cA_|)$ in these spaces. 
\ \ $\Box$

\BEL
\label{lem6.3}
An  orthonormal basis in $A_{V,\eta}^2(\cA_|)$ is formed by the functions
\bea
f_{n,\eta} (z) = C_{n,\eta} \frac{z^{n + \eta/(2\pi)}}{v(z)}, \ \ n \in \bbZ,  \label{6.20}
\eea
where $C_{n,\eta}$ are the normalization constants,
\bea
C_{n,\eta}^{-2} &= & \int\limits_\cA \Big| \frac{z^{n + \eta/(2\pi)}}{v(z)} 
\Big|^2 V(z)\, dA(z) = \int\limits_0^{2 \pi} 
\int\limits_{1/\rho}^\rho r^{2n + \eta/\pi +1 } dr d \theta 
\roweq
\frac{2 \pi}{2(n+1)  + \eta/\pi  } ( \rho^{2(n+1) + \eta/\pi } 
- \rho^{-2(n+1)-\eta/\pi })  . \label{6.22}
\eea
\ENL

\bigskip

Proof. The mutual orthogonality of the functions \eqref{6.20} follows from the 
usual orthogonality  relation of the monomials $z^n$ and the fact that the factor  
$1/v$ is cancelled by the weight in the inner product of $A_V^2(\cA)$. Indeed,  
given $n,m \in \bbZ$ and $\eta \in [-\pi,\pi]$ we have
\bea
& &\int\limits_\cA f_{n,\eta} \overline{f_{m,\eta}} V dA
= C_{n,\eta} C_{m,\eta} \int\limits_\cA z^{n  + \eta/(2 \pi) } \overline{z^{m +
\eta/(2\pi)}}dA   \label{6.89}
\eea
where $ z^{\eta/(2 \pi) } = (r e^{i \theta})^{\eta/(2 \pi) }= r^{ \eta/(2 \pi) } 
e^{i \theta (  \eta/(2 \pi)  + \eta \ell (\theta) ) }$ for some 
$\ell (\theta) \in \bbZ$ hence, the integral in \eqref{6.89} reads in 
polar coordinates as
\beas
\int\limits_{1/\rho}^\rho r^{\ldots} dr \int\limits_0^{2 \pi} 
e^{i \theta ( n - m + \eta/(2 \pi)  + \eta \ell (\theta)- \eta/(2 \pi)- 
\eta \ell (\theta) ) } d \theta
\eeas
which is null unless $n=m$.

As for the completeness of the orthonormal sequence \eqref{6.20}, 
an arbitrary $h \in A^2(\cA)$ can be approximated
in the Bergman space $A^2(\cA)$ by a linear combination of the functions $z^n$, 
$n\in \bbZ$, since these form an orthonormal basis of $A^2(\cA)$, after a proper 
normalization. Since the mapping $f \mapsto  v^{-1} f$ is an isometry from 
$A^2(\cA)$ onto $A_V^2(\cA)$, we find that an arbitrary $g \in A_V^2(\cA)$ can be 
approximated in $A_V^2(\cA)$ by a linear combination of the functions $v^{-1}z^n$, 
$n\in \bbZ$.
In view of \eqref{3.27b}, the linear combinations of the functions \eqref{6.20}
are dense in $A_{V,\eta, {\rm ext}}^2(\cA)$ and thus the system
\eqref{6.20} is complete in $A_{V,\eta}^2(\cA_|)$.  \ \ $\Box$

\bigskip

We now construct the kernel function $K_\eta(z,w)$, see \eqref{1.17}. Due to Lemma 
\ref{lem6.3}, the kernel $K_{\eta, \cA}$ of the orthogonal projection 
$f \mapsto \int_\cA K_{\eta, \cA} ( \cdot, w) f(w) dA(w)$ from $L_V^2(\cA)$ onto
$\cA_{V,\eta}^2(\cA_|)$ is 
\bea
K_{\eta, \cA} (z,w)= \sum_{n \in \bbZ} f_{n, \eta} (z) 
\overline{f_{n, \eta} (w) } V(w),
\eea
where the weight function $V$ comes from the inner product of $L_V^2(\cA)$.  
We have, by \eqref{3.24r}, $I_2^{-1} I_1^{-1} f  = f \circ \phi \circ E $  
hence, Lemma \ref{lem2.5} and \eqref{3.24r} imply that  the kernel 
$K_\eta$ of the orthogonal  
projection from $L^2(\varpi)$ onto $A_\eta^2(\varpi)$ can be written as
\bea
& & K_\eta(z,w) = K_{\eta, \cA} \big( \phi\circ E(z), \phi \circ E(w) \big) 
\big| (\phi\circ E)' (w) \big|^2 
\roweq 
\sum_{n \in \bbZ} \frac{C_{n, \eta}^2
\phi( e^{ i 2 \pi z} )^{n + \eta/(2\pi)} 
\overline{\phi( e^{ i 2 \pi  w } )}^{n + \eta/(2\pi)}}{
v( \phi( e^{ i 2 \pi z} ) ) \overline{ v( \phi( e^{ i 2 \pi w} ) )} } ,
\ \ \ \  z,w \in \varpi,
\eea
where we used \eqref{3.24}. Moreover, we have
\beas
v( \phi( e^{ i 2 \pi z} ) )  
= \frac{1}{2 \pi}  \frac{e^{ - i 2 \pi z}}{\phi'(e^{ i 2 \pi z}) } 
\eeas
so that denoting
\bea
\widetilde K(z,w) &=& 
4 \pi^2 e^{ i 2 \pi (z-  \bar w)} 
\phi'( e^{ i 2 \pi z} )\overline{ \phi'( e^{  i 2 \pi  w} )  }
\label{6.23}
\eea
the kernel becomes
\bea
\ \ K_\eta(z,w) = 
\widetilde K(z,w) 
\sum_{n \in \bbZ} C_{n, \eta}^2 \phi( e^{ i 2 \pi z} )^{n + \eta/(2\pi)} 
\overline{\phi( e^{ i 2 \pi \bar w } )}^{n + \eta/(2\pi)} ,
\ \   z,w \in \varpi . 
\label{6.24}
\eea

Next, we recall from Lemma \ref{lem0.2} the conformal mapping 
\bea
\varphi(z)   
= \frac{1}{i 2 \pi} \log \big( \phi ( e^{i 2 \pi z} ) \big) + [{\rm Re} \, z]
 \ \ \ \ \mbox{with}  \ \ \phi ( e^{i 2 \pi z} ) = e^{ i 2 \pi \varphi(z) }. 
\label{6.26}
\eea
Taking into account \eqref{6.22} we now  write \eqref{6.24} as
\bea
 K_\eta(z,w) &=& \widetilde K(z,w)
\sum_{n \in \bbZ}\frac{2n + \eta/\pi }{ 2\pi( \rho^{2n + \eta/\pi  } - 
\rho^{-2 n - \eta/\pi } ) }  e^{ i (2 \pi (n-1) + \eta)( 
\varphi(z) -  \overline{ \varphi(w) } )  }.  
\label{6.30}
\eea
By \eqref{6.26}, we have  $|e^{ i 2 \pi \varphi(z)}| = 
|\phi( e^{i 2 \pi z})|$ for all $z \in  \varpi$, and since 
$\phi(e^{i 2 \pi z}) \in \cA$, we have
\bea
1 /\rho < \big| e^{ i 2 \pi \varphi(z)} \big| < \rho ,
\eea
which implies that the series \eqref{6.30} converges absolutely and uniformly on
all compact subsets of $\varpi$. We arrive at the main result of this section.

\BET \label{th3.8}
Let the domain  $\varpi$ (equivalently, $\Pi$) be  simply connected and 
satisfy assumption $(A)$ of Section \ref{sec0}, and let $\eta \in [-\pi, \pi]$.  The  
kernel $K_\eta$ of the projection from $L^2(\varpi)$ onto $A_\eta^2(\varpi)$, see 
\eqref{1.17}, is given by formulas \eqref{6.23}, \eqref{6.30}. 
\ENT

\section{Applications.}
\label{sec5}

In order to demonstrate the potential applicability of Theorem \ref{th3.8}, we
derive the kernel formula \eqref{0.6} from it. To this end
we combine the kernel formula \eqref{6.30} with  \eqref{4.4zz} in Theorem 
\ref{th2.3}, 
\bea
& & K_\Pi (z,w) 
=  \frac{ \widetilde K(z,w)}{2 \pi} 
\int\limits_{-\pi}^{\pi}  
\sum_{n \in \bbZ} 
\frac{2n + \eta/\pi }{ 2\pi( \rho^{2n + \eta/\pi  } - \rho^{-2n - \eta/\pi} ) }  
\nonumber 
\\ & & \times 
e^{ i \eta \big( [{\rm Re}\,z]- [{\rm Re}\,w] \big)
+ i (2 \pi (n-1)  + \eta) \big( 
\varphi(z - [{\rm Re}\,z]) - \varphi( \bar w - [{\rm Re}\,w])   \big)}  d\eta
\roweq   
\frac{ \widetilde K(z,w)}{4 \pi^2}
 e^{-i 2 \pi( \varphi(z) -  \overline{ \varphi(w) } ))  } 
4\pi \int\limits_{-\infty}^\infty  \frac{ t}{ \rho^{2t} - 
\rho^{ - 2 t } } 
e^{ i 2 \pi t   ( \varphi(z) - \overline{  \varphi(w) )} } d t  .
\label{3.41b}
\eea
where we also made the summation and integration over $[-\pi, \pi]$ into an 
integration over the real line (with $t = \eta /(2 \pi)$) and used  $\varphi\big(z - [{\rm Re} \, z] \big) 
=\varphi(z)- [{\rm Re} \, z]$ (see \eqref{6.26}) and similarly for the variable 
$w$. There holds the integral formula (Fourier transform)
\bea
\ \ \ \int\limits_{- \infty}^\infty \frac{t}{2} \, {\rm csch}\,(at) e^{-i st} dt = 
\int\limits_{- \infty}^\infty  \frac{t}{e^{at} - e^{-at} } e^{-i st} dt = 
\frac{1}{4a^2} \pi^2 {\rm sech}^2 \Big( \frac{\pi s}{2a} \Big) , \  s \in \bbR, \label{3.84}
\eea
where $a>0$ is a parameter and csch denotes the hyperbolic cosecant; the identity \eqref{3.84} can be got 
from the  Fourier cosine transform given in the table 1.9(18) of \cite{tables}, 
since the transformed function is even. We apply \eqref{3.84} with
$a = 2  \log \rho $ and 
$ 
s=  2\pi( \varphi(z) -  \overline{ \varphi( w)} ) 
$. 
The equations 
\eqref{6.23}, \eqref{3.41b} yield formula \eqref{0.6}:

\BEC
\label{th3.10} If the periodic domain $\Pi$, \eqref{1.2}, is simply connected
and satisfies the assumption $(A)$ of Section \ref{sec0},
then its Bergman  kernel equals 
\bea
\ \ \ \ K_\Pi (z,w)  &=& \widetilde K_\Pi(z,w) \frac{\pi^3}{ 4 (\log \rho)^2 } 
{\rm sech}^2 \Big( \frac{ \pi^2( \varphi(z) - \overline{ \varphi( w)} ) }{
2 \log \rho}  \Big)  ,  
\label{3.50}
\eea
where 
\bea
\widetilde K_\Pi(z,w) =
e^{ i 2 \pi (z-   \varphi(z) - \bar w +   \overline{ \varphi(w) } ) } 
\phi'( e^{ i 2 \pi z} )\overline{ \phi'( e^{  i 2 \pi  w} )  },
 \label{3.50s}
\eea
$\phi: \cD \to \cA$ is the Riemann mapping between doubly 
connected domains $\cD = \{ e^{i 2 \pi z}\, : \, z \in \varpi\}$ and the 
annulus $\cA = \{ w \, : \,  \rho^{-1 } < w < \rho \}$, $\rho > 1$, and 
$\varphi: \Pi \to \bbC$ is given in \eqref{6.26}.
\ENC

Note that in the case $\Pi$ is the strip 
$\Sigma = (-\infty,\infty) \times (-\pi,\pi)$, see the proof of Proposition
\ref{prop0.4}, the mapping $\phi$ 
is the identity, and obviously formula \eqref{3.50}--\eqref{3.50s}, with 
log\,$\rho = 2 \pi^2$,  boils down to \eqref{BS}. 

The invariance of the kernel under the mapping
\bea
(z,w) \mapsto (z + m , w+ m)  , \ \ \ m \in \bbZ , \label{3.50f}
\eea
follows from \eqref{6.26}. 
We also have the following estimate for the kernel; estimates for arbitrary 
$(z,w) \in \Pi \times \Pi$ follow by applying the mapping \eqref{3.50f}.

\BEC
\label{cor6.1} Let $\Pi$ be as in Corollary \ref{th3.10} and assume in addition 
that there exists a constant $C >0$ such that 
\bea
\frac1C \leq \sup_{w \in \cD} |\phi'(w)| \leq C.  \label{3.50t}
\eea
There exist constants $0 < C_1 < C_2$ such that 
\bea
C_1 e^{- \pi^2 n/(2 \log \rho)} \leq |K_\Pi(z,0)| \leq 
C_2 e^{- \pi^2 n/(2 \log \rho)}
\label{3.50u}
\eea
for all $n \in \bbN$ and all $z \in \Pi$ with $ n -1\leq |z|  \leq n$.
\ENC

Proof. We have for all $z \in \Pi$
\beas
\exp\Big(\frac{\pi^2 \varphi(z)}{2 \log \rho } \Big)= 
\exp\Big(\frac{\pi^2  \varphi(z - [{\rm Re}\, z ])}{2 \log \rho } \Big)
\exp\Big(\frac{\pi^2  [{\rm Re}\, z ]}{2 \log \rho } \Big),
\eeas
and on the right hand side  the modulus of the first factor is bounded
and bounded from below by a positive constant, since $\varphi$ is a maps
$\varpi$ onto a bounded domain. The same holds for the factor 
$\widetilde K_\Pi$, by
the assumption \eqref{3.50t}. The result follows from the properties
of sech. \ \ $\Box$

\BER
Condition \eqref{3.50t} holds, if the boundary of the doubly connected
domain $\cD$, equivalently, the boundary of $\Pi$,  is for example $C^2$-smooth.  
There are non-smooth domains $\cD$ such that for some boundary point $\zeta \in 
\partial \cD$ we have $\lim_{z \to \zeta} |\phi'(\zeta)| = \infty$; see the
references mentioned below \eqref{3.20}. In this case the decay \eqref{3.50u} 
of course fails, since \eqref{3.50},  \eqref{3.50s} imply
$
\sup_{z \in \varpi_m} |K_\Pi(z,0)| = \infty $ for all
$m \in \bbZ. $ 
\ENR

The following decay estimate for the kernel does not depend on the 
geometry of the periodic cell $\varpi$. 

\BEC
\label{cor6.2} Let $\Pi$ be as in Corollary \ref{th3.10} and let $G \subset 
\varpi$ be a compact subset. Then, there exist constants $0 < C_1 < C_2$ 
such that 
\bea
C_1 e^{- \pi^2  |n|/(2 \log \rho)} \leq |K(z,0)| \leq  
C_2 e^{- \pi^2  |n|/(2 \log \rho)}
\eea
for all $n \in \bbZ$ and all $z \in \Pi$ such that $z - n \in G$.
\ENC

\BER
If the periodic cell $\varpi$ were doubly connected, the entire domain
$\Pi$ would be infinitely connected, and in principle the above method 
yields a representation for the Bergman kernel $K_\Pi$ in terms
of Riemann mappings between finitely instead of infinitely connected domains. 
However, 
obtaining the  concrete formula \eqref{3.50}--\eqref{3.50s} required
the use of a quite explicit, simple sequence of orthogonal functions
on the annulus and an exact calculation of their norms. The author does
not know a good enough example of such sequences for triply and higher 
connected domains. 
\ENR

We finally prove a simple boundedness result for the Bergman projection
with respect to certain weighted
$L^p$-norms. Let us consider continuous weights $W: \Pi \to \bbR^+$ which 
only depend on the real part of the variable $z \in \Pi$.
We assume that there are constants $a, C >0$ and $0< b<1$
such that for all $x \in \bbR$, $n \in \bbZ$
\bea
\frac1C W(x) e^{ - a |n|^b } \leq |W(x + n ) | \leq C W(x) e^{ a |n|^b }  . \label{7.2}
\eea

\BET
\label{th5.2}
Let $\Pi$ be as in Corollary \ref{cor6.1}.
Let  $W: \Pi \to \bbR^+$ be a weight as above and $1 \leq p < \infty$.
Then, the projection operator $P_\Pi: L_W^p(\Pi) \to L_W^p(\Pi)$ is bounded.
\ENT

Proof. 
We apply the Schur test, \cite{Z}, Theorem 3.6.,
with the constant test function $h(z) \equiv 1$; note that we
take $K_\Pi(z,w) W(w)^{-1}$ for the kernel $K$ in the reference.
By also taking into account that $|\widetilde K_\Pi|$ is function bounded 
by a constant $C > 0$ in $\Pi\times \Pi$,
we obtain from \eqref{3.50u}, \eqref{7.2} for all $z = x + iy \in \Pi$
\bea
& & \int\limits_{ \Pi} 
| K_\Pi(z,w) | W(w)^{-1} W(z)   dA(z) 
\rowleq C \sum_{n \in \bbN} \int\limits_{n-1 \leq |x - \xi| \leq n} 
\exp \Big( - \frac{ \pi^2 n  }{2 \log \rho} 
\Big) \exp \big( a n^b) W(x)^{-1} W(x)    dA(w) 
\rowleq 
C'  \sum_{n \in \bbN} 
\exp \Big( - \frac{ \pi^2 n  }{2 \log \rho} + a n^b \Big) 
\rowleq C'' .  \label{7.10}
\eea
The second integral $ \int_{ \Pi} 
| K_\Pi(z,w) |  dA(w) $ in the Schur test is also bounded by a 
constant; the estimation is easier. 
\ \ $\Box$

\vfill \eject

\end{document}